\documentclass{amsart}
\usepackage{amscd,amssymb,relsize, hyperref,verbatim}
\usepackage[arrow,matrix,graph,frame,poly,arc,tips]{xy}

\newtheorem{thm}{Theorem}[section]
\newtheorem*{theo}{Theorem}
\newtheorem{cor}[thm]{Corollary}

\newtheorem{prop}[thm]{Proposition}

\theoremstyle{definition}
\newtheorem{definition}[thm]{Definition}
\newtheorem{example}[thm]{Example}
\newtheorem{remark}[thm]{Remark}

\numberwithin{equation}{section}

\newcommand{\abs}[1]{\left|#1\right|}
\newcommand{\set}[1]{\left\{#1\right\}}

\renewcommand{\b}[1]{\mathbf{#1}}

\newcommand{\A}{\mathcal{A}}

\newcommand{\cV}{\mathcal{V}}
\newcommand{\cE}{\mathcal{E}}
\newcommand{\cC}{\mathcal{C}}

\newcommand{\bS}{\mathbb{S}}
\newcommand{\bT}{\mathbb{T}}
\newcommand{\R}{\mathbb{R}}
\newcommand{\C}{\mathbb{C}}
\newcommand{\Q}{\mathbb{Q}}

\newcommand{\Z}{\mathbb{Z}}
\renewcommand{\k}{\Bbbk}

\DeclareMathOperator{\id}{id}

\DeclareMathOperator{\tc}{{\sf TC}}
\DeclareMathOperator{\zcl}{\sf zcl}
\DeclareMathOperator{\cat}{cat}
\DeclareMathOperator{\secat}{secat}
\DeclareMathOperator{\cl}{\sf cl}

\begin{document}

\title[Motion planning in tori]%
{Motion planning in tori}

\author[D.~C. Cohen]{Daniel C. Cohen$^\dag$}
\address{Department of Mathematics, Louisiana State University,
Baton Rouge, LA 70803}
\email{\href{mailto:cohen@math.lsu.edu}{cohen@math.lsu.edu}}
\urladdr{\href{http://www.math.lsu.edu/~cohen/}
{http://www.math.lsu.edu/\~{}cohen}}
\thanks{{$^\dag$}Partially supported 
by National Security Agency grant H98230-05-1-0055}

\author[G. Pruidze]{Goderdzi Pruidze}
\address{Department of Mathematics, Louisiana State University,
Baton Rouge, LA 70803}
\email{\href{mailto:gio@math.lsu.edu}{gio@math.lsu.edu}}
\urladdr{\href{http://www.math.lsu.edu/~gio/}%
{http://www.math.lsu.edu/\~{}gio}}



\begin{abstract}
Let $X$ be a subcomplex of the standard CW-decomposition of the $n$-dimensional torus.  We exhibit an explicit optimal motion planning algorithm for $X$.  This construction is used to calculate the topological complexity of complements of general position arrangements and Eilenberg-Mac\,Lane spaces associated to right-angled Artin groups.
\end{abstract}

\date{April 27, 2007}

\maketitle

\section{Introduction} \label{sec:intro}

Let $X$ be a path-connected topological space.  We assume throughout that $X$ has the homotopy type of a finite CW-complex.  Viewing $X$ as the space of configurations of a mechanical system, the motion planning problem consists of constructing an algorithm which takes as input pairs of configurations $(x_0,x_1) \in X \times X$, and produces a continuous path $\gamma\colon [0,1] \to X$ from the initial configuration $x_0=\gamma(0)$ to the terminal configuration $x_1=\gamma(1)$.  The motion planning problem is of interest in robotics, see, for example, Latombe \cite{La} and Sharir \cite{Sh}.

In a recent sequence of papers, Farber \cite{Fa03,Fa04,Fa05} develops a topological approach to the motion planning problem.  Let $PX$ be the space of all continuous paths $\gamma\colon [0,1] \to X$, equipped with the compact-open topology.  The map $\pi\colon PX \to X \times X$ 
defined by sending a path to its endpoints, $\pi\colon\gamma\mapsto (\gamma(0),\gamma(1))$, is a fibration, with fiber $\Omega{X}$, the based loop space of $X$.  The motion planning problem then asks for a section of this fibration, a map $s\colon X\times X \to PX$ satisfying $\pi \circ s = \id_{X\times X}$.

It would be desirable for the motion planning algorithm to depend continuously on the input.  However, it is not difficult to show that there exists a globally continuous motion planning $s\colon X\times X \to PX$ if and only if $X$ is contractible, see \cite[Thm.~1]{Fa03}.  This leads to the following definition, see \cite[\S{2}]{Fa04}.

\begin{definition} \label{def:motion planner}
A \emph{motion planner} for $X$ is a collection of subsets $F_1,\dots,F_k$ of $X\times X$ and continuous maps $s_i\colon F_i \to PX$ such that
\begin{enumerate}
\item the sets $F_i$ are pairwise disjoint, $F_i \cap F_j=\emptyset$ if $i \neq j$, and cover $X \times X$, 
\[
X\times X = F_1 \cup \cdots \cup F_k;
\]
\item $\pi \circ s_i = \id_{F_i}$ for each $i$; and
\item each $F_i$ is a Euclidean Neighborhood Retract.
\end{enumerate}
\end{definition}
\noindent The sets $F_i$ are referred to as the \emph{local domains} of the motion planner, and the maps $s_i$ are the \emph{local rules}.  

Call a motion planner optimal if it requires a minimal number of local domains (resp., rules).  
Define the \emph{topological complexity} of $X$, $\tc(X)$, to be the number of local domains in an optimal motion planner for $X$.  In this paper, we determine the topological complexity of an arbitrary  subcomplex of the standard CW-decomposition of the $n$-dimensional torus $(S^1)^{\times n}$.  Examples include the skeleta of the $n$-torus, as well as cell complexes associated to right-angled Artin groups.  We calculate the topological complexity of these spaces by explicitly exhibiting optimal motion planners.  Our methods also yield analogous motion planners for standard subcomplexes of $\bS^{\times n}$, where $\bS$ is a sphere of any dimension.

In \cite[\S{31}]{Fa05}, Farber poses the problem of computing the topological complexity of an aspherical space $X$, a space for which the higher homotopy groups vanish, $\pi_i(X)=0$ for $i\ge 2$.  We solve this problem in the case where $X=K(G,1)$ is an Eilenberg-Mac\,Lane space associated to an arbitrary right-angled Artin group $G$.  

Let $\Gamma=(\cV_\Gamma,\cE_\Gamma)$ be a finite graph with vertex set $\cV_\Gamma$, edge set $\cE_\Gamma \subseteq\binom{\cV_\Gamma}{2}$, and no loops or multiple edges.  The \emph{right-angled Artin group}, or graph group, associated to $\Gamma$ is the group $G=G_\Gamma$ with generators $v$ corresponding to vertices $v \in \cV_\Gamma$ of $\Gamma$, and commutator relators $vw=wv$ corresponding to edges $\set{v,w}\in \cE_\Gamma$ of $\Gamma$.  Let $X_\Gamma$ be the subcomplex of the standard CW-decomposition of the $n$-dimensional torus $(S^1)^{\times n}$ obtained by deleting the cells corresponding to noncliques of $\Gamma$.  As shown by Charney and Davis \cite{CD} and Meier and Van Wyk \cite{MVW}, the complex $X_\Gamma$ is an Eilenberg-Mac\,Lane space of type $K(G_\Gamma,1)$.

\begin{theo}
Let $\Gamma$ be a finite simplicial graph.  Then the topological complexity of the associated cell complex $X_\Gamma$ is $\tc(X_\Gamma)=z(\Gamma)+1$, where 
\[
z(\Gamma) = \max_{K_1,K_2} \left| \cV_{K_1} \cup \cV_{K_2} \right|
\]
is the largest number of vertices of $\Gamma$ covered by precisely two cliques $K_1$ and $K_2$ 
in $\Gamma$.
\end{theo}

Let $\A$ be an arrangement of hyperplanes in $\C^\ell$, a finite collection of codimension one affine subspaces.  A principal object in the topological study of arrangements is the \emph{complement} 
$X_\A=\C^\ell \setminus \bigcup_{H\in \A}H$.  The study of the topological complexity of complements of arrangements was initiated by Farber and Yuzvinsky \cite{FY}, who determined $\tc(X_\A)$ for numerous arrangements, including the braid arrangement, with complement $X_\A=F_\ell(\C)$, the configuration space of $\ell$ ordered points in $\C$.  Their results imply that, for an arrangement $\A$ of $n$ hyperplanes in general position in $\C^\ell$, if $n \ge 2\ell$, then $\tc(X_\A)=2\ell+1$.  Left open by these results is the case of an arbitrary general position arrangement.

By a classical result of Hattori \cite{Ha}, the complement $X_\A$ of an arrangement $\A$ of $n> \ell$ general position hyperplanes in $\C^\ell$ has the homotopy type of the $\ell$-dimensional skeleton 
$(S^1)^{\times n}_\ell$ 
of the $n$-dimensional torus.  (If $n\le\ell$, then $X_\A \simeq (S^1)^{\times n}= (S^1)^{\times n}_n$ has the homotopy type of the $n$-torus.)  Since topological complexity is a homotopy-type invariant (see \cite[Thm. 3]{Fa03} and Section \ref{sec:tc} below), we have $\tc(X_\A)=\tc( (S^1)^{\times n}_\ell)$.

\begin{theo}
Let $\A$ be an arrangement of $n$ hyperplanes in general position in $\C^\ell$. Then the topological complexity of the complement $X_\A$ is $\tc(X_\A)=\min\set{n+1,2\ell+1}$.
\end{theo}

This recovers a recent result of Yuzvinsky \cite{Yuz07}.

The cell complexes $X_\Gamma$ associated to right-angled Artin groups and complements $X_\A$  of (general position) hyperplane arrangements are known to be formal spaces, in the sense of Sullivan.  Our results support the conjecture (stated explicitly for arrangements in \cite{Yuz07}) that the topological complexity of a formal space $X$ is determined by the cohomology ring of $X$.

This note is organized as follows.  After briefly recalling some relevant results on topological complexity and motion planning in Section \ref{sec:tc}, we produce an optimal motion planner for an arbitrary subcomplex of the standard CW-decomposition of $\bS^{\times n}$, where $\bS$ is an odd-dimensional sphere, in Section \ref{sec:tori}.  Cell complexes associated to right-angled Artin groups and complements 
of general position hyperplane arrangements, as well as higher-dimensional analogues, are analyzed in Section \ref{sec:ra} and Section \ref{sec:gp} respectively.  In Section \ref{sec:even}, we exhibit an optimal motion planner for an arbitrary subcomplex of the standard CW-decomposition of $\b{S}^{\times n}$, where $\b{S}$ is an even-dimensional sphere.

\section{Topological complexity and motion planning} \label{sec:tc}

Let $p\colon E\to B$ be a fibration.  The 
\emph{sectional category}, or Schwarz genus, of $p$, denoted by $\secat(p)$, 
is the smallest integer $k$ such that $B$ can be covered 
by $k$ open subsets, over each of which $p$ has a continuous section.  
See James \cite{Ja} as a classical reference.  
The topological nature of the motion planning problem is made clear by the following result.

\begin{thm}[{\cite[cf.~\S{13}]{Fa05}}] \label{thm:tc=secat} 
Let $X$ be a simplicial polyhedron.  Then the topological complexity of $X$ is equal to the sectional category of the path-space fibration $\pi\colon PX \to X \times X$,  
\\[5pt]
\centerline{\phantom{\qed} \hfill 
$\tc(X)=\secat(\pi)$. \hfill \qed
}
\end{thm}

In \cite[Thm.~3]{Fa03}, Farber shows that the sectional category of the path-space fibration 
$\pi\colon PX \to X \times X$, and hence the topological complexity of $X$, is an invariant of the homotopy type of $X$.  From the realization $\tc(X)=\secat(\pi)$, one also obtains a number of bounds on the topological complexity of $X$ in terms of the dimension and the Lusternik-Schnirelman category $\cat(X)$.  For instance, one has
\[
\cat(X) \le \tc(X) \le 2 \cat(X)-1 \le 2 \dim(X)+1.
\]
We will not have cause to use these bounds.  
We will, however, make extensive use of a cohomological lower bound provided by the realization 
$\tc(X)=\secat(\pi)$.

Let $\k$ be a field, and let $A=\bigoplus_{k=0}^\ell A^k$ be a graded $\k$-algebra with $A^k$ finite-dimensional for each $k$.  Define the cup length of $A$, denoted by $\cl(A)$, to be the largest integer $q$ for which there are homogeneous elements $a_1,\dots,a_q$ of positive degree in $A$ such that $a_1\cdots a_q \neq 0$.  If $p\colon E \to B$ is a fibration, the sectional category admits the following lower bound:
\[
\secat(p) > \cl\bigl(\ker(p^*\colon H^*(B;\k) \to H^*(E;\k)\bigr),
\]
see \cite[\S{8}]{Ja}.  We subsequently work exclusively with a fixed field $\k$ of characteristic zero, and suppress coefficients in cohomology, writing $H^*(Y)=H^*(Y;\k)$ for ease of notation.

For the path-space fibration $\pi\colon PX \to X\times X$, using the fact that $PX\simeq X$ and the 
K\"unneth formula, the kernel of $\pi^*\colon H^*(X\times X)
\to H^*(PX)$ may be identified with the kernel $Z=Z(H^*(X))$ 
of the cup-product map $H^*(X)\otimes H^*(X)
\xrightarrow{\ \cup\,} H^*(X)$, see \cite[Thm.~7]{Fa03}.  
Call the cup length of the ideal $Z$ of zero divisors
the \emph{zero-divisor cup length} of $H^*(X)$, and write $\zcl(H^*(X))=\cl(Z)$.  

\begin{prop} \label{prop:zcl}
The topological complexity of $X$ is greater than the zero-divisor cup length of $H^*(X)$, 
\\[5pt]
\centerline{\phantom{\qed} \hfill 
$\tc(X)> \zcl(H^*(X))$. \hfill \qed
}
\end{prop}

We conclude this section by recalling the construction of a motion planner for a product space $X\times Y$ from \cite[\S{12}]{Fa04}.  Suppose that $X$ admits a motion planner $X\times X=F_1\cup\dots\cup F_k$, $f_i\colon F_i \to PX$, for which $F_1\cup \dots \cup F_i$ is closed for each $i$, $1\le i \le k$.  Analogously, suppose $Y$ admits a motion planner $Y\times Y=G_1\cup\dots\cup G_\ell$, $g_j\colon G_j \to PY$, with $G_1\cup \dots \cup G_j$ closed for each $j$, $1\le j \le \ell$.  
Identify $(X\times Y)\times (X\times Y)=X\times X\times Y \times Y$.  
For each $r$, $2 \le r \le k+\ell$, let 
\[
V_r=\bigcup_{i+j=r} F_i \times G_j,
\]
and define $h_r\colon V_r \to P(X\times Y)$ by $h_r|_{F_i \times G_j} = f_i \times g_j$.

\begin{prop} \label{prop:product}
The subsets $V_r\subset (X\times Y) \times (X\times Y)$ and maps $h_r\colon V_r \to P(X\times Y)$, $2 \le r \le k+\ell$, are the local domains and rules of a motion planner for $X\times Y$.  Hence, 
\\[6pt]
\centerline{\phantom{\qed} \hfill 
$\tc(X\times Y)\le \tc(X) + \tc(Y) - 1$. \hfill \qed
}
\end{prop}

\section{Motion planning in subcomplexes of tori} \label{sec:tori}

Let $\bS$ be an odd-dimensional sphere.  View $\bS=S^{2k-1}$ as the set of all points $x \in \C^k$
 with $|x|=1$.  Analogously, view $\bS^{\times n}=\bS\times\dots\times\bS$ as the set of all points $\b{x}=(x_1,\dots,x_n)\in (\C^k)^n$ such that $|x_i|=1$ for each $i$.  Let $e=(1,0,\dots,0)$ be the (first) standard unit vector in $\C^k$.  
We will work with the standard CW-decomposition of $\bS^{\times n}$, with cells given by 
\begin{equation} \label{eq:sphere cells}
\bS^{\times n}_J=\{\b{x}\in \bS^{\times n} \mid x_i=e \text{ if } j\notin J, x_i \neq e
\text{ if } j \in J \}
\end{equation}
for subsets $J$ of $[n]=\set{1,\dots,n}$.  

The purpose of this section is to exhibit an optimal motion planner for an arbitrary subcomplex $X$ of this CW-decomposition of $\bS^{\times n}$.  Subcomplexes of products of even-dimensional spheres behave differently, and are treated in Section \ref{sec:even}.  
While the primary focus in the applications presented in Sections \ref{sec:ra} and \ref{sec:gp} 
will be on subcomplexes of the $n$-dimensional torus $(S^1)^{\times n}$, the construction is independent of (odd) dimension.  So we present the general case.  We begin by recalling a motion planner for $\bS$ itself.
 
\begin{example}[{\cite[Ex.~7.5]{Fa05}}] \label{ex:odd sphere}
The sphere $\bS=S^{2k-1}$ admits a motion planner with $2$ local domains.  Let $F_1 \subset \bS\times \bS$ be the subset consisting of all pairs of antipodal points, $F_1=\{(x,-x) \mid x \in \bS\}$, and define $s_1\colon F_1 \to P\bS$ as follows:  Fix the standard nowhere zero tangent vector field $\nu$ on $\bS$, and move $x$ towards $-x$ along the semicircle tangent to the vector $\nu(x)$.

For a second local domain, let $F_2=\bS\times \bS \setminus F_1$ be the complement of $F_1$ in $\bS\times \bS$.  For $(x,y) \in F_2$, we have $x \neq -y$ and we may define $s_2\colon F_2 \to P\bS$ by moving $x$ towards $y$ along the shortest geodesic arc.

A calculation in the cohomology ring reveals that $\zcl(H^*(\bS))=1$, so $\tc(\bS)=2$ and the above motion planner is optimal.
\end{example}

\begin{remark} \label{rem:even sphere}
The topological complexity of an even-dimensional sphere is $\tc(S^{2k})=3$, see \cite[Thm.~8]{Fa03}. 
An optimal motion planner for $S^{2k}$ is exhibited in Example \ref{ex:even sphere}.
\end{remark}

Note that the set $F_1$ of pairs of antipodal points is closed in $\bS\times \bS$.  Consequently, we may apply the algorithm of \cite[\S12]{Fa04} to obtain a motion planner for the product space $\bS^{\times 2}=\bS\times\bS$, see Proposition \ref{prop:product} and the preceding discussion.  Iterating this construction, we obtain a motion planner for 
$\bS^{\times n}$.   For a subset $I \subseteq[n]$, define 
\begin{equation} \label{eq:sub domain}
F_I = \{(\b{x},\b{y}) \in \bS^{\times n}  \times \bS^{\times n}  \mid x_i=-y_i \text{ if } i \in I, x_i\neq -y_i \text{ if } i \notin I\},
\end{equation}
and define $s_I\colon F_I \to P\bS^{\times n} $ using the maps $s_1$ and $s_2$ given in Example \ref{ex:odd sphere} coordinate-wise.  For $(\b{x},\b{y})=\bigl((x_1,y_1),\dots,(x_n,y_n)\bigr)$, set
\[
s_I(\b{x},\b{y}) = \bigl(t_1(x_1,y_1),\dots,t_n(x_n,y_n)\bigr) \in P\bS \times\dots\times P\bS = 
P\bS^{\times n} ,
\]
where $t_i = s_1$ if $i \in I$, and $t_i=s_2$ if $i \notin I$.  

For $j=0,1,\dots,n$, let
\begin{equation*} \label{eq:torus domains}
W_j = \bigcup_{|I|=n-j}F_I,
\end{equation*}
and define $\sigma_j\colon W_j \to P\bS^{\times n}$ by $\sigma_j|_{F_I} = s_I$.

\begin{prop} \label{prop:torus motion plan}
The subsets $W_j \subset \bS^{\times n} \times \bS^{\times n}$ and maps 
$\sigma_j\colon W_j \to P\bS^{\times n}$ are the local domains and rules of an 
optimal motion planner for 
$\bS^{\times n}$.  Hence, $\tc(\bS^{\times n})=n+1$.
\end{prop}
\begin{proof}
It is clear from the above construction that the subsets $W_j \subset \bS^{\times n} \times \bS^{\times n}$ and maps $\sigma_j\colon W_j \to P\bS^{\times n}$ comprise a motion planner for $\bS^{\times n}$, so $\tc(\bS^{\times n}) \le n+1$.

If $\bS$ is the $(2k-1)$-sphere, the cohomology ring $H^*(\bS^{\times n})$ is an exterior algebra on degree $2k-1$ generators $e_1,\dots,e_n$.  Let $\bar{e}_i=1\otimes e_i-e_i\otimes 1 \in
H^*(\bS) \otimes H^*(\bS)$, and note that $\bar{e}_i$ is a zero-divisor.  Checking that $\prod_{i=1}^n \bar{e}_i \neq 0$ in $H^*(\bS^{\times n}\times\bS^{\times n})$, we have $\zcl(H^*(\bS^{\times n}))\ge n$.  Hence, as shown in \cite[Thm.~13]{Fa03}, the topological complexity of $\bS^{\times n}$ is equal to $n+1$.  Consequently, the motion planner constructed above is optimal.
\end{proof}

Now let $X$ be a subcomplex of the standard CW-decomposition of $\bS^{\times n}$.  
We will use the motion planner for $\bS^{\times n}$ constructed above to produce an optimal motion planner for $X$.  First, we need a definition.   
For $J \subset [n]$, let $\bT_J \cong \bS^{\times |J|}$ denote the subcomplex  of 
$\bS^{\times n}$ given by 
$\bT_J = \set{  \b{x} \in \bS^{\times n} \mid x_i = e\text{ if } i \notin J}$.  
If $J$ and $K$ are disjoint subsets of $[n]$, then $\bT_J \cap \bT_K=(e,\dots,e)$ is 
a point, and $\bT_J \cup \bT_K = \bT_J \vee \bT_K$ is the wedge of $\bT_J$ and $\bT_K$, 
and is a subcomplex of $\bS^{\times n}$.
Define
\begin{equation*} \label{eq:z(X)}
z(X)=\max\set{ \abs{J}+\abs{K} \mid J \cap K=\emptyset
\text{ and }\bT_J \vee \bT_K\text{ is a  subcomplex of } X}.
\end{equation*}
Note that, for instance, the set $K$ may be empty.

\begin{thm} \label{thm:subcomplex}
Let $X$ be a subcomplex of the standard CW-decomposition of $\bS^{\times n}$.  
Then $X$ admits an optimal motion planner with $z(X)+1$ local domains.  Hence, 
$\tc(X)=z(X)+1$.
\end{thm}

We will establish this result by (i) exhibiting a motion planner for $X$ with $z(X)+1$ local domains, and (ii) showing that the zero-divisor cup length of $H^*(X)$ is at least $z(X)$.

\begin{prop} \label{prop:subcomplex}
Let $X$ be a subcomplex of the standard CW-decomposition of $\bS^{\times n}$.  
Then $X$ admits a motion planner with at most $z(X)+1$ local domains.
\end{prop}
\begin{proof} 
Let $\sigma_j(X)$ denote the restriction of the map $\sigma_j\colon W_j \to P\bS^{\times n}$ to 
$(X\times X) \cap W_j$.  
We assert that the decomposition $X \times X = \bigcup_{j=n-z(X)}^n (X \times X) \cap W_j$, 
together with the maps $\sigma_j(X)$, $n-z(X)\le j\le n$, comprise a motion planner for $X$.  
The cells of $X$ are of the form $A_J = \bS^{\times n}_J$ for certain subsets $J$ of $[n]$.  Consequently, $X\times X$ admits a CW-decomposition with cells $A_J \times A_K$ for certain subsets $J$ and $K$ of $[n]$.  Note that if $A_J \times A_K$ is a cell of $X \times X$, then $\abs{J\cup K} \le z(X)$.

Consider the intersection of such a cell $A_J \times A_K$ with one of the subsets $F_I$ defined 
in \eqref{eq:sub domain} above.  If there exists $i \in I$ with $i \notin J \cup K$, it is readily checked that $(A_J \times A_K)\cap F_I=\emptyset$.  Explicitly, since $i \in I$, we have $x_i=-y_i$ in $F_I$.  Since $i\notin J$ and $i\notin K$, we have $x_i=e$ and $y_i=e$ in $A_J \times A_K$.  Consequently, $(A_J \times A_K)\cap F_I=\emptyset$ 
if $I\not\subset J \cup K$.

If $I \subset J \cup K$, then $(A_J \times A_K)\cap F_I$ consists of those points $(\b{x},\b{y})$ which satisfy $x_i=y_i=e$ for $i \notin J \cup K$, and the following conditons for $i \in J \cup K$:
\[
\begin{matrix}
x_i=-y_i, x_i\neq e, y_i\neq e\ \text{if $i \in J$, $i\in K $, $i\in I$,} \hfill &
x_i\neq -y_i, x_i\neq e, y_i\neq e\ \text{if $i \in J\cap K$, $i\notin I$,} \hfill  \\
x_i=-e, y_i=e \  \text{if $i \in J$, $i\notin K$, $i\in I$,}\hfill\   &
x_i\neq -e,y_i=e \ \text{if $i \in J$, $i\notin K$, $i\notin I$,}\hfill \\
x_i=e, y_i=-e \  \text{if $i \notin J$, $i \in K$, $i\in I$,} \hfill &
x_i=e, y_i\neq -e \ \text{if $i \notin J$, $i \in K$,  $i\notin I$.} \hfill
\end{matrix}
\]
In particular, $(A_J \times A_K)\cap F_I\neq \emptyset$ if $I\subset J\cup K$.

The above observations may be used to show that $(X \times X) \cap W_j = \emptyset$ if $j < n-z(X)$.  
Fix such a $j$, and recall that $W_j=\bigcup_{|I|=n-j} F_I$.  If $A_J \times A_K$ is a cell of $X \times X$, then $|J \cup K| \le z(X) <  |I|$ for each $I$ with $|I|=n-j$.  Consequently, for each such $I$, there is an $i \in I$ with $i \notin J$ and $i\notin K$, and $(A_J \times A_K) \cap F_I =\emptyset$.
It follows that $(X \times X) \cap W_j = \emptyset$ for $j < n-z(X)$.  

If $\abs{I} \ge n-z(X)$, let $\phi_I$ denote the restriction of the map $s_I\colon F_I \to P\bS^{\times n}$ 
to $(A_J \times A_K)\cap F_I$.  
For $(\b{x},\b{y}) \in \bS^{\times n} \times \bS^{\times n}$, $\phi_I(\b{x},\b{y})$ is a particular path taking $\b{x}$ to $\b{y}$, that is, $x_i$ to $y_i$ for each $i$.  Recall that $x_i$ is taken to $y_i$ in the $i$-th coordinate sphere along the geodesic determined by the (fixed) vector field $\nu$ if $x_i= -y_i$, and along the shortest geodesic if $x_i\neq -y_i$.  Note that the latter is simply a constant path if $x_i=y_i$.  The above description of $(A_J \times A_K)\cap F_I$ may be used to check that the image of 
$\phi_I$ is contained in $PX$.  If $i \in J$ (resp., $i \in K$), then the $i$-th coordinate sphere is contained in the closure of the cell $A_J$ (resp., $A_K$) in $X$.  If $i \notin J$ and $i\notin K$, then $x_i=y_i=e$, and the $i$-th component of the path $\phi_I(\b{x},\b{y})$ is constant.  Consequently, if $(\b{x},\b{y})\in X \times X$, then $\phi_I(\b{x},\b{y})$ is a path from $\b{x}$ to $\b{y}$ in $X$.  
It follows that 
the sets $(X\times X) \cap W_j$ and maps $\sigma_j(X)\colon (X\times X) \cap W_j \to PX$, $n-z(X) \le j \le n$,  comprise a motion planner for $X$ as asserted.
\end{proof}

The above result implies that $\tc(X) \le z(X)+1$.  We establish the reverse inequality by analyzing the cohomology ring $H^*(X)$.  As noted in the proof of Proposition \ref{prop:torus motion plan}, 
if $\bS=S^{2k-1}$, the cohomology ring $H^*(\bS^{\times n})$ is an exterior algebra $E$ on degree $2k-1$ generators $e_1,\dots,e_n$.  For $J=\set{j_1,\dots,j_k} \subset [n]$, let $e_J = e_{j_1}\cdots e_{j_k} \in E$, and recall that $\bS_J^{\times n}$ denotes the corresponding cell of $\bS^{\times n}$, see \eqref{eq:sphere cells}.  The following result was established for cell complexes associated to right-angled Artin groups by Charney and Davis \cite[Thm.~3.2.4]{CD}.  Their argument may be brought to bear on an arbitrary subcomplex of a product of odd-dimensional spheres.

\begin{prop} \label{prop:H*X}
If $X$ is a subcomplex of the standard CW-decomposition of $\bS^{\times n}$, then the cohomology ring $H^*(X)$ is the quotient of the exterior algebra $E$ by the ideal
\\[6pt]
\centerline{\phantom{\qed} \hfill 
$I_X=\langle e_J \mid J \subset [n]\text{ and } \bS_J^{\times n}\text{ is not a cell of } X\rangle$. \hfill \qed
}
\end{prop}

\begin{prop} \label{prop:zclX}
If $X$ is a subcomplex of the standard CW-decomposition of $\bS^{\times n}$, then the zero-divisor cup length of $H^*(X)$ is at least $z(X)$.
\end{prop}
\begin{proof}
Write $z=z(X)$ and assume that $J_1=\set{1,\dots,k}$ and $J_2=\set{k+1,\dots, z}$ are the underlying index sets of the subcomplex $\bT_{J_1} \vee \bT_{J_2}$ of $X$ realizing $z(X)=z$.  First, consider the element $\bar{e}_1\bar{e}_2\cdots\bar{e}_z$ in $E\otimes E$, where $E$ is the exterior algebra and $\bar{e}_i=1\otimes e_i - e_i\otimes 1$.  As noted in \cite[Lemma 10]{FY} in greater generality (when $E$ is generated in degree $1$), one can check that
\begin{equation} \label{eq:prod1}
\bar{e}_1\bar{e}_2\cdots\bar{e}_z = \sum_{(J,J')} (-1)^{|J|}\mathrm{sign}(\sigma)\, e_J \otimes e_{J'},
\end{equation}
where the sum is over all partitions $(J,J')$ by ordered subsets of $[z]=\set{1,\dots,z}$, and $\sigma$ is the shuffle on $[z]$ which puts every element of $J'$ after all elements of $J$, preserving the orders inside $J$ and $J'$.

Now consider the cohomology ring $H^*(X)=E/I_X$ of the subcomplex $X$. We assert that the (image of the) element $\bar{e}_1\bar{e}_2\cdots\bar{e}_z$ is non-zero in $H^*(X) \otimes H^*(X) = E/I_X \otimes E/I_X$.  For this, first note that $\pm e_{J_1}\otimes e_{J_2}=\pm e_1\cdots e_k \otimes e_{k+1}\cdots e_z$ is a summand of \eqref{eq:prod1}, and that $e_{J_1}$ and $e_{J_2}$ are non-zero in $H^*(X)$ since the epimorphic images are non-zero in the cohomology rings $H^*(\bT_{J_1})$ and $H^*(\bT_{J_2})$ respectively.  Next, observe that the projection $E\otimes E \to E/I_X \otimes E/I_X$ is given by annihilating multiples of elements of the form $e_K\otimes 1$ and $1\otimes e_K$ if $\bS^{\times n}_K$ is not a cell in $X$.  Consequently, the image of the sum \eqref{eq:prod1} in $E/I_X \otimes E/I_X$ is given by
\[
\sum_{(J,J')} (-1)^{|J|}\mathrm{sign}(\sigma)\, e_J \otimes e_{J'}
\]
where the sum is over all partitions $(J,J')$ as above for which there are corresponding subcomplexes 
$\bT_J$  and $\bT_{J'}$ of $X$.  In particular, the summand $\pm e_{J_1} \otimes e_{J_2}$ survives, and the zero-divisor cup length of $H^*(X)= E/I_X$ is at least $z=z(X)$.
\end{proof}

Together, Proposition \ref{prop:subcomplex} and Proposition \ref{prop:zclX} imply that the topological complexity of a subcomplex $X$ of the standard CW-decomposition of $\bS^{\times n}$ is 
$\tc(X)=z(X)+1$, thereby proving Theorem \ref{thm:subcomplex}.  We note two consequences.

\begin{cor} \label{cor:optimal zcl}
Let $X$ be a subcomplex of the standard CW-decomposition of $\bS^{\times n}$.  Then 
\begin{enumerate}
\item the motion planner $\sigma_j(X)\colon (X\times X) \cap W_j \to PX$, $n-z(X) \le j \le n$, 
is optimal, and
\item the zero-divisor cup length of the cohomology ring $H^*(X)$ is $\zcl(H^*(X))=z(X)$. \qed
\end{enumerate}
\end{cor}

We conclude this section with two brief applications.

In general, the topological complexity of a product space $X \times Y$ admits the upper bound 
$\tc(X\times Y) \le \tc(X)+\tc(Y)-1$, see Proposition \ref{prop:product}.  
In the case where $X$ and $Y$ are subcomplexes of $\bS^{\times n}$ equality always holds.

\begin{prop} \label{prop:x product}
If $X_1$ and $X_2$ are subcomplexes of $\bS^{\times n_1}$ and $\bS^{\times n_2}$, then
\[
\tc(X_1\times X_2)=
\tc(X_1)+\tc(X_2)-1.
\]
\end{prop}
\begin{proof} 
By Theorem \ref{thm:subcomplex}, it suffices to show that $z(X_1 \times X_2)\ge z(X_1)+z(X_2)$. 
For $i=1,2$, let $\bT_{J_i} \vee \bT_{K_i}$ be a subcomplex of $X_i$ realizing $z(X_i)=\abs{J_i}+\abs{K_i}$.  Then the product space $X_1\times X_2$ contains $(\bT_{J_1} \times \bT_{J_2}) \vee (\bT_{K_1} \times \bT_{K_2}) = \bT_{J}\vee\bT_{K}$ as a subcomplex, where $J\subset [n_1+n_2]$ is the subset corresponding to $J_1\subset[n_1]$ and $J_2\subset[n_2]$ in the obvious manner, and similarly for $K$.  Thus, 
\[
z(X_1 \times X_2) \ge \abs{J}+\abs{K}= (\abs{J_1}+\abs{J_2}) +(\abs{K_1}+\abs{K_2}) = z(X_1)+z(X_2)
\]
as required. 
\end{proof}

A similar result holds for wedges.  In general, the topological complexity of a wedge $X \vee Y$ admits the following upper bound
\begin{equation*} \label{eq:xwedge}
\tc(X \vee Y ) \le \max\set{\tc(X),\tc(Y),\cat(X)+\cat(Y)-1},
\end{equation*}
see \cite[Thm.~19.1]{Fa05}.  In the case where $X$ and $Y$ are subcomplexes of $\bS^{\times n}$, equality always holds.

\begin{prop} \label{prop:x wedge}
If $X_1$ and $X_2$ are subcomplexes of $\bS^{\times n_1}$ and $\bS^{\times n_2}$, then
\[
\tc(X_1\vee X_2)=
\max\set{\tc(X_1),\tc(X_2),\cat(X_1)+\cat(X_2)-1}.
\]
\end{prop}
\begin{proof}
For a subcomplex $X$ of $\bS^{\times n}$, define
\[
d(X) = \max\set{\abs{J} \mid \bT_J\text{ is a subcomplex of } X}.
\]
Note that if $\bS=S^{2k-1}$, then $\dim(X) = (2k-1)\cdot d$.  Note also that $X$ is $(2k-2)$-connected.

For the wedge $X=X_1 \vee X_2$, observe that $z(X) \ge d(X_1)+d(X_2)$.  
If $z(X)$ is realized by a subcomplex $\bT_{J_1} \vee \bT_{K_1}$ of $X_1$, then $z(X)=z(X_1)$, and similarly for $X_2$.  If neither of these possibilities occurs, then $z(X)$ is realized by a wedge 
$\bT_J\vee\bT_K$, where $\bT_J$ is a subcomplex of $X_1$ and $\bT_K$ is a subcomplex of $X_2$.
Hence, $z(X)=\max\set{z(X_1),z(X_2),d(X_1)+d(X_2)}$.  
By Theorem \ref{thm:subcomplex}, we have
\[
\tc(X)=\tc(X_1\vee X_2)=
\max\set{z(X_1)+1,z(X_2)+1,d(X_1)+d(X_2)+1}.
\]

Checking that the cup-length in $H^*(X_i)$ is equal to $d(X_i)$, we have $\cat(X_i) = d(X_i)+1$ using the classical inequalities $\cl(H^*(X))+1\le\cat(X) \le 1+\dim(X)/r$ for an $(r-1)$-connected space $X$, see James \cite{Ja}.  Thus, 
$\cat(X_1)+\cat(X_2)-1=d(X_1)+d(X_2)+1$. 
The result follows.
\end{proof}

\section{Right-angled Artin groups} \label{sec:ra}
Let $\Gamma$ be a finite simplicial graph, with vertex set $\cV$ and edge set $\cE\subset \binom{\cV}{2}$.  Associated to $\Gamma$ is a right-angled Artin group $G_\Gamma$, with a generator $v$ corresponding to each vertex $v \in\cV$, and a commutator relation $vw=wv$ for each edge $\{v,w\}\in\cE$.  That is, the group $G_\Gamma$ has presentation
\[
G_\Gamma = \langle v \in \cV \mid vw=wv \text{ if } \{v,w\}\in\cE\rangle.
\]

Let $\Delta_\Gamma$ be the flag complex associated to the graph $\Gamma$.  The flag complex is the largest simplicial complex with $1$-skeleton equal to $\Gamma$.  The $k$-simplices of $\Delta_\Gamma$ are the $(k+1)$-cliques of $\Gamma$, i.e., the complete subgraphs on $k+1$ vertices in $\Gamma$.

Let $n=|\cV|$ be the cardinality of the vertex set of $\Gamma$, and consider the $n$-torus $(S^1)^{\times n}$ with its standard CW-decomposition.  
Denote by $X_\Gamma$ 
the subcomplex of $(S^1)^{\times n}$ obtained by deleting the cells corresponding to the non-faces of the flag complex $\Delta_\Gamma$.  The 
complex $X_\Gamma$ is an Eilenberg-Mac\,Lane space of type $K(G_\Gamma,1)$, see \cite{CD,MVW}.

Let $c_k(\Gamma)$ denote the number of $k$-cliques in $\Gamma$ (where $c_0(\Gamma)=1$).  It is readily checked that the integral homology group $H_k(X_\Gamma;\Z)$ is free abelian of rank $c_k(\Gamma)$.  Consequently, the Poincar\'e polynomial of $X_\Gamma$ 
is equal to the clique polynomial of the graph
\[
P(X_\Gamma,t)=\sum_{k\ge 0} b_k(X_\Gamma)t^k=\sum_{k\ge 0} c_k(\Gamma)t^k=C(\Gamma,t).
\]
If the cardinality of the largest clique in $\Gamma$ is equal to $d$, then $\dim X_\Gamma=d$.

Label the vertices of $\Gamma$, and write $\cV=\{v_1,\dots,v_n\}$.  
The cohomology ring $H^*(X_\Gamma)$ 
is the exterior Stanley-Reisner ring of $\Gamma$, see \cite{CD,KR}.  Explicitly, $H^*(X_\Gamma)=E/I_\Gamma$ is the quotient of the exterior algebra $E$ on generators $e_1,\dots,e_n$ in degree one by the ideal 
$I_\Gamma=\langle e_i e_j \mid \{v_i,v_j\}\notin \cE\rangle$.

Denote by $z(\Gamma)$ the maximal number of vertices of $\Gamma$ that are covered by precisely two cliques.  That is, $z(\Gamma)$ is the largest cardinality of union of the vertex sets of two cliques in $\Gamma$.  
Observe that $z(\Gamma)$ may be realized by disjoint cliques in $\Gamma$.

\begin{thm} \label{thm:right-angled tc}
If $\Gamma$ is a finite simplicial graph, then the topological complexity of the associated cell complex $X_\Gamma$ is $\tc(X_\Gamma)=z(\Gamma)+1$.
\end{thm}
\begin{proof}
By Theorem \ref{thm:subcomplex}, it suffices to show that $z(\Gamma)=z(X_\Gamma)$.  
Identify the vertex set $\cV$ of $\Gamma$ with $[n]=\set{1,\dots,n}$, and let $J,K\subset [n]$ be disjoint cliques realizing $z(\Gamma)$.  Then $J$ and $K$ correspond to faces of the flag complex $\Delta_\Gamma$, and it follows that $\bT_J\vee\bT_K$ is a subcomplex of $X_\Gamma$.  Consequently, $z(\Gamma) \le z(X_\Gamma)$.

If $\bT_J\vee\bT_K$ is a subcomplex of $X_\Gamma$, then $J$ and $K$ are disjoint, so correspond to disjoint cliques in $\Gamma$.  This implies that $z(X_\Gamma)\le z(\Gamma)$, so we have $z(\Gamma)=z(X_\Gamma)$ as required.
\end{proof}

\begin{cor} \label{cor:right-angled zcl}
Let $\Gamma$ be a finite simplicial graph with associated cell complex $X_\Gamma$.  Then the zero-divisor cup length of the cohomology ring $H^*(X_\Gamma)$ is $\zcl(H^*(X_\Gamma))=z(\Gamma)$.
\qed
\end{cor}

The cell complex $X_\Gamma$ associated to the graph $\Gamma$ has higher-dimensional analogues $X_\Gamma^k$ introduced by Kim and Roush \cite{KR}.  Fix a positive integer $k$, and recall that $\bS=S^{2k-1}$ is the $(2k-1)$-sphere.  Given a simplicial graph $\Gamma$ with $|\cV|=n$, let $X^k_\Gamma$ be the subcomplex of the standard CW-decomposition of $\bS^{\times n}$ obtained by deleting the cells corresponding to the non-faces of the flag complex $\Delta_\Gamma$.  Note that $X_\Gamma^1=X_\Gamma$ is the $K(G_\Gamma,1)$-space considered previously, and that $X_\Gamma^k$ is simply-connected for $k \ge 2$.

The cohomology ring $H^*(X^k_\Gamma)$ is a rescaling of $H^*(X_\Gamma)=E/I_\Gamma$ in the sense of \cite{PS04}.  Explicitly, if $E^k$ denotes the exterior algebra generated by elements $e_1^k,\dots,e_n^k$ of degree $2k-1$, then $H^*(X^k_\Gamma)=E^k/I^k_\Gamma$, where $I^k_\Gamma$ is the ideal generated by $\set{e^k_i e^k_j \mid \set{v_i,v_j}\notin\cE}$.  It follows that the zero-divisor cup length of $H^*(X^k_\Gamma)$ is equal to that of $H^*(X_\Gamma)$, 
\[
\zcl(H^*(X^k_\Gamma))=\zcl(H^*(X_\Gamma))=z(\Gamma).
\]
Furthermore, since $X^k_\Gamma$ is a subcomplex of the standard CW-decomposition of $\bS^{\times n}$, Theorem \ref{thm:subcomplex} provides a motion planner for $X^k_\Gamma$ with local domains $(X^k_\Gamma \times X^k_\Gamma) \cap W_j$, $n-z(X_\Gamma^k) \le j \le n$.  Clearly, $z(X_\Gamma^k)=z(\Gamma)$.  
Hence, we have the following generalization of Theorem \ref{thm:right-angled tc}.

\begin{thm} \label{thm:high dim ra tc}
For any positive integer $k$, the topological complexity of the complex $X^k_\Gamma$ is equal to $z(\Gamma)+1$. 
\qed
\end{thm}

\begin{remark} \label{rem:ra formal}
The complexes $X^k_\Gamma$ are known to be formal spaces, see \cite{PS06}.  Since $X^k_\Gamma$ is simply-connected for $k\ge 2$, results of \cite{FGKV,LM} imply that the topological complexity of the rationalization $(X^k_\Gamma)_\Q$ of $X^k_\Gamma$ is also equal to $z(\Gamma)+1$ 
in this instance.
\end{remark}

\section{General position arrangements} \label{sec:gp}

Let $\A$ be an arrangement of $n$ affine hyperplanes in general position in $\C^\ell$, with complement 
$X_\A=\C^\ell \setminus \bigcup_{H \in \A} H$.
If $n\le \ell$, it is readily checked that the complement has the homotopy type of the $n$-dimensional torus, $X_\A \simeq (S^1)^{\times n}$.  
If $n>\ell$, then by a classical result of Hattori \cite{Ha}, the complement has the homotopy type of the $\ell$-dimensional skeleton of the $n$-torus, $X_\A \simeq (S^1)^{\times n}_\ell$.  
For $n \ge 2\ell$, results of Farber and Yuzvinsky \cite{FY} imply that the topological complexity of $X_\A$ is equal to $2\ell+1$.  For arbitrary $n$, we have the following.

\begin{thm} \label{thm:tc gp} 
Let $\A$ be an arrangement of $n$ hyperplanes in general position in $\C^\ell$. Then the topological complexity of the complement of $\A$ is $\tc(X_\A)=\min\set{n+1,2\ell+1}$.
\end{thm}
\begin{proof} If $n \le \ell$, then the topological complexity of $X_\A \simeq (S^1)^{\times n}$ is $n+1$, see Proposition \ref{prop:torus motion plan}.

If $n> \ell$, since $X_\A\simeq (S^1)^{\times n}_\ell$ by Hattori's theorem, we have 
$\tc(X_\A)=\tc((S^1)^{\times n}_\ell)=z((S^1)^{\times n}_\ell)+1$ by Theorem \ref{thm:subcomplex}.  
So it suffices to show that $z((S^1)^{\times n}_\ell)=\min\set{n,2\ell}$.

Let $z=\min\set{n,2\ell}$, and consider the subsets $J_0=\set{1,\dots,\ell}$ and $K_0=\set{\ell+1,\dots,z}$ 
of $[n]=\set{1,\dots,n}$.  Then the $\ell$-skeleton $(S^1)^{\times n}_\ell$ of the $n$-torus contains the wedge of tori $\bT_{J_0}\vee\bT_{K_0}$ as a subcomplex, so  $z((S^1)^{\times n}_\ell)\ge z=\min\set{n,2\ell}$.

To establish the reverse inequality, assume 
that $(S^1)^{\times n}_\ell$ contains a wedge of tori $\bT_J\vee\bT_K$ as a subcomplex, 
where $J$ and $K$ are disjoint subsets of $[n]$ and  
$\abs{J}+\abs{K}>z=\min\set{n,2\ell}$.  Since $(S^1)^{\times n}_\ell$ is $\ell$-dimensional, the subsets $J$ and $K$ of $[n]$ are of cardinality at most $\ell$.  So we have 
$\min\set{n,2\ell}<\abs{J}+\abs{K}\le 2\ell$, and $\min\set{n,2\ell}=n$.  Since $\abs{J}+\abs{K}>n$, the subsets $J$ and $K$ of $[n]$ cannot be disjoint.  This contradiction completes the proof.
\end{proof}

\begin{cor} \label{cor:zcl gp} 
Let $\A$ be an arrangement of $n$  general position hyperplanes in $\C^\ell$ with complement 
$X_\A$. Then the zero-divisor cup length of the cohomology ring $H^*(X_\A)$ is 
$\zcl(H^*(X_\A))=\min\set{n,2\ell}$. \qed
\end{cor}

A hyperplane arrangement $\cC$ in $\C^\ell$ is said to be central if every hyperplane of $\cC$ contains the origin.  A central arrangement $\cC$ is generic if every $\ell$ element subset of the hyperplanes of $\cC$ is independent.  Note that $|\cC| \ge \ell$ for $\cC$ generic. In other words, a generic arrangement $\cC$ of $n$ hyperplanes in $\C^\ell$ is the cone of an affine arrangement $\A$ of $n-1$ general position hyperplanes in $\C^{\ell-1}$.  See Orlik and Terao \cite{OT} as a general reference on arrangements.

\begin{thm} \label{thm:tc gen}
Let $\cC$ be a generic arrangement of $n$ hyperplanes in $\C^\ell$. Then the topological complexity of the complement of $\cC$ is $\tc(X_\cC)=\min\set{n+1,2\ell}$.
\end{thm}
\begin{proof}
As noted above, the generic arrangement $\cC\subset\C^\ell$ is the cone of an affine arrangement $\A$ of $n-1$ general position hyperplanes in $\C^{\ell-1}$.  The relationship between the complements is well known:  We have $X_\cC \cong X_\A \times \C^*$, see \cite[Prop.~5.1]{OT}.   Since $X_\A \simeq (S^1)^{\times n-1}_{\ell-1}$ has the homotopy type of the $(\ell-1)$-skeleton of the $(n-1)$-torus and $\C^*\simeq S^1$, it follows that $X_\cC \simeq (S^1)^{\times n-1}_{\ell-1} \times S^1$.  Using Theorem \ref{thm:tc gp} and Proposition \ref{prop:x product}, we obtain
\[
\tc(X_\cC) = \tc((S^1)^{\times n-1}_{\ell-1})+\tc(S^1)-1 = \tc(X_\A)+1=\min\set{n,2\ell-1}+1=\min\set{n+1,2\ell}
\]
as asserted. 
\end{proof}

\begin{cor} \label{cor:zcl gen} 
Let $\cC$ be a generic arrangement of $n$ hyperplanes in $\C^\ell$ with complement 
$X_\A$. Then the zero-divisor cup length of the cohomology ring $H^*(X_\cC)$ is 
$\zcl(H^*(X_\cC))=\min\set{n,2\ell-1}$. \qed
\end{cor}

\begin{remark} \label{rem:yuz}
After proving Theorems \ref{thm:tc gp} and \ref{thm:tc gen}, we learned that these results had been established previously by Yuzvinsky \cite{Yuz07}.  The arguments are similar.
\end{remark}

\begin{remark} \label{rem:nilpotent formal}
In \cite{LM}, Lechuga and Murillo show that the topological complexity of the rationalization $X_\Q$ of a nilpotent formal space $X$ is given by $\tc(X_\Q)=\zcl(H^*(X_\Q))+1$.  
It is well known that complements of complex hyperplane arrangements are formal.  If $\A$ is an 
arrangement of $n$ general position hyperplanes in $\C^\ell$ with $\ell \ge 2$ (resp., $\cC$ is a 
generic arrangement in $\C^\ell$ with $\ell\ge 3$), the fundamental group $G$ of the complement is free 
abelian, hence is nilpotent. However, the higher homotopy groups are not, in general, nilpotent as $\Z{G}$-modules, see \cite{PS02}.
\end{remark}

Following \cite{CCX}, to an arrangement $\A$ of hyperplanes in $\C^\ell$, we associate a 
redundant arrangement $\A^{k}$ of codimension $k$ subspaces in $\C^{k\ell}$.  
Given a hyperplane
$H\subset\C^\ell$, let $H^{k}$ be the codimension $k$ affine subspace
of $\C^{k\ell}=(\C^{\ell})^{k}$ consisting of all $k$-tuples of points
in $\C^{\ell}$, each of which lies in $H$.  Applying this construction to each of 
the hyperplanes in 
$\A$ yields an arrangement $\A^k$ of 
codimension $k$ subspaces in $\C^{k\ell}$, with complement
$X_\A^k=\C^{k\ell}\setminus \bigcup_{H\in\A} H^k$.

The cohomology ring $H^*(X_\A^k)$ of the complement of the redundant subspace arrangement is a rescaling of the cohomology ring $H^*(X_\A)$ of the original hyperplane arrangement.  By the Orlik-Solomon theorem, $H^*(X_\A)=E/I_\A$ is the quotient of an exterior algebra on $n=|\A|$ degree one generators $e_1,\dots,e_n$ by a certain homogeneous ideal $I_\A$, see \cite[\S\S{3.1--3.2}]{OT}.  If $E^k$ is the exterior algebra generated by elements $e_1^k,\dots,e_n^k$ of degree $2k-1$, and $I_\A^k$ is the ideal in $E^k$ corresponding to the ideal $I_\A$ of $E$ in the obvious manner, then $H^*(X_\A^k)=E^k/I_\A^k$.  See \cite[\S{2}]{CCX} for details.  It follows that the zero-divisor cup length of $H^*(X_\A^k)$ is equal to that of $H^*(X_\A)$, $\zcl(H^*(X_\A^k))=\zcl(H^*(X_\A))$ for each $k$.

\begin{thm} \label{thm:tc gp red} 
Let $\A$ be an arrangement of $n$ hyperplanes in general position in $\C^\ell$. Then, for each $k$, the topological complexity of the complement of the redundant subspace arrangement $\A^k$ is $\tc(X_\A^k)=\min\set{n+1,2\ell+1}$.
\end{thm}
\begin{proof}
For $k=1$, this is Theorem \ref{thm:tc gp}.  So assume that $k \ge 2$.  Write $\bS=S^{2k-1}$.

If $n \le \ell$, checking that $X^k_\A \simeq \bS^{\times n}$, we have $\tc(X^k_\A)=n+1$, see Proposition \ref{prop:torus motion plan}.

For $n>\ell$, as shown by Papadima and Suciu \cite[Lemma 8.4]{PS04}, the analogue of Hattori's theorem holds for redundant subspace arrangements associated to general position arrangements.  If $\A$ is a general position arrangement of $n$ hyperplanes in $\C^\ell$, then the complement of the subspace arrangement $\A^k$ has the homotopy type of the $(2k-1)\ell$-skeleton 
$(\bS^{\times n})_\ell$ of the standard CW-decomposition of $\bS^{\times n}$.  Since
$X_\A^k \simeq (\bS^{\times n})_\ell$, by Theorem \ref{thm:subcomplex}, it suffices to show that 
$z((\bS^{\times n})_\ell)=\min\set{n,2\ell}$.  This may be accomplished by generalizing the argument given in the proof of Theorem \ref{thm:tc gp} in the obvious manner.
\end{proof}

If $u$ and $v$ are two distinct points in a manifold $X$, the \emph{open string configuration space} is
\begin{equation*} \label{eq:open string}
G_n(X,u,v)=\set{(x_1,\dots,x_n) \in X^{\times n} \mid x_1 \neq u,\ x_n \neq v,\ x_i \neq x_{i+1}
\ \text{for}\ 1\le i\le n-1}.
\end{equation*}
These configuration spaces arise naturally in the topological study of billiard problems, see~\cite{Fa02}.

\begin{prop} \label{prop:open string}
The open string configuration space $G_n(\C^k,u,v)$ is equal to the complement $X_\A^k$ of a redundant subspace arrangement associated to an arrangement $\A$ of $n+1$ hyperplanes in general position in $\C^n$.
\end{prop}
\begin{proof}
Choose coordinates $y_1,\dots,y_n$ on $\C^n$, and let $\A\subset \C^n$ be the arrangement of $n+1$ hyperplanes defined by the polynomial $Q(\A)=y_1(y_n-1)\prod_{i=1}^{n-1}(y_i-y_{i+1})$.  It is readily checked that the hyperplanes of $\A$ are in general position.

Without loss of generality, take $u=(0,\dots,0)$ and $v=(1,\dots,1)$ in $\C^k$.  
It is then a straightforward exercise to check that $G_n(\C^k,u,v)$ coincides with the complement of the subspace arrangement $\A^k$ associated to the hyperplane arrangement $\A$ in $\C^n$ defined above.
\end{proof}

This result, together with Theorem \ref{thm:tc gp red}, yields the following.

\begin{cor} \label{cor:open string tc}
The topological complexity of the open string configuration space $G_n(\C^k,u,v)$ is equal to $n+2$. \qed
\end{cor}

\section{Products of even-dimensional spheres} \label{sec:even}
In this section, we consider the case of an even-dimensional sphere $\b{S}=S^{2k}$.  View the product space $\b{S}^{\times n}=\b{S}\times\dots\times\b{S}$ as the set of points $\b{x}=(x_1,\dots,x_n)\in(\R^{2k+1})^n$ such that $\abs{x_i}=1$ for each $i$, and let $e=(1,0,\dots,0)\in\R^{2k+1}$.  Let $X$ be a subcomplex of the standard CW-decomposition of $\b{S}^{\times n}$, the latter with cells 
$\b{S}^{\times n}_J=\set{\b{x} \in \b{S}^{\times n} \mid x_i=e \text{ if } j\notin J, x_i \neq e
\text{ if } j \in J}$ for subsets $J$ of $[n]$.  We will construct an optimal motion planner for $X$.

\begin{example}[{\cite[Ex.~7.5]{Fa05}}] \label{ex:even sphere}
The sphere $\b{S}=S^{2k}$ admits a motion planner with $3$ local domains. Let $F_0=\set{(e,-e)}\subset\b{S}\times\b{S}$, where $e$ is the fixed point given above, and choose the local rule $s_0$ to be any fixed path from $e$ to $-e$. 

Let $F_1 \subset \b{S}\times \b{S}$ be the subset consisting of all pairs of antipodal points except 
$(e, -e)$,  $ F_1=\set{(x,-x) \in \b{S}\times \b{S}\mid x\ne e} $, and define $s_1\colon F_1 \to P\b{S}$ as follows.   Fix a nowhere zero tangent vector field $\nu$ on $\b{S} \setminus \set{e}$, and move $x$ towards $-x$ along the semicircle tangent to the vector $\nu(x)$.

For a third local domain, let $F_2=\b{S}\times \b{S} \setminus (F_0 \cup F_1) =\set{(x,y) \in \b{S}  \times \b{S} \mid  x\ne -y}$ be the complement of the union of $F_0$ and $F_1$ in $\b{S}\times \b{S}$.  For $(x,y) \in F_2$, we have $x \neq -y$ and we may define $s_2\colon F_2 \to P\b{S}$ by moving $x$ towards $y$ along the shortest geodesic arc.

One can check that $\zcl(H^*(\b{S}))=2$, so $\tc(\b{S})=3$ and the above motion planner is optimal.
\end{example}

Applying the algorithm of \cite[\S{12}]{Fa04} summarized in Proposition \ref{prop:product}, 
we obtain a motion planner for $\b{S}^{\times n}$.  Since we now have choices of $3$ local domains for each sphere $\b{S}$, we parameterize (subsets of) local domains of $\b{S}^{\times n}$ using multi-indices $\alpha=(\alpha(1),\dots,\alpha(n))$, where $\alpha(i) \in\set{0,1,2}$, that is, using functions 
$\alpha\colon [n] \to \set{0,1,2}$.  Given $\alpha$, define 
\[
F_{\alpha} = \set{(\b{x},\b{y}) \in \b{S}^{\times n}  \times \b{S}^{\times n}  \mid (x_i,y_i)\in F_{\alpha(i)} },
\]
and define $s_{\alpha}\colon F_{\alpha} \to P\b{S}^{\times n} $ using the maps $s_0$, $s_1$ and $s_2$ given in Example \ref{ex:even sphere} coordinate-wise.  For $(\b{x},\b{y})=\bigl((x_1,y_1),\dots,(x_n,y_n)\bigr)$, set
\[
s_{\alpha}(\b{x},\b{y}) = \bigl(t_1(x_1,y_1),\dots,t_n(x_n,y_n)\bigr) \in P\b{S} \times\dots\times P\b{S} = 
P\b{S}^{\times n} ,
\]
where $t_i = s_{\alpha (i)}$.  

For $j=0,1,\dots,2n$, let
\begin{equation*} \label{eq:even torus domains}
W_j = \bigcup_{|\alpha|=j}F_{\alpha},
\end{equation*}
where $\abs{\alpha} = \alpha(1) + \dots + \alpha(n)$, and define $\sigma_j\colon W_j \to P\b{S}^{\times n}$ by $\sigma_j|_{F_{\alpha}} = s_{\alpha}$.

\begin{prop} \label{prop:even torus motion plan}
The subsets $W_j \subset \b{S}^{\times n} \times \b{S}^{\times n}$ and maps 
$\sigma_j\colon W_j \to P\b{S}^{\times n}$ are the local domains and rules of an 
optimal motion planner for 
$\b{S}^{\times n}$.  Hence, $\tc(\b{S}^{\times n})=2n+1$.
\end{prop}
\begin{proof}
It is clear from the above construction that the subsets $W_j \subset \b{S}^{\times n} \times \b{S}^{\times n}$ and maps $\sigma_j\colon W_j \to P\b{S}^{\times n}$ comprise a motion planner for $\b{S}^{\times n}$, so $\tc(\b{S}^{\times n}) \le 2n+1$.

If $\b{S}$ is the $2k$-sphere, the cohomology ring $H^*(\b{S}^{\times n})$ is a commutative algebra of square-free monomials on degree $2k$ generators $e_1,\dots,e_n$.  
The zero-divisor $\bar{e}_i=1\otimes e_i-e_i\otimes 1$ in $H^*(\b{S}) \otimes H^*(\b{S})$ 
satisfies $(\bar{e}_i)^2=-2\bar{e}_i\otimes\bar{e}_i \neq 0$ since we work with coefficients in a field $\k$ of characteristic zero.  
Checking that $\prod_{i=1}^n (\bar{e}_i)^2 \neq 0$ in $H^*(\b{S}^{\times n}\times\b{S}^{\times n})$, we have $\zcl(H^*(\b{S}^{\times n}))\ge 2n$.  Hence, as shown in \cite[Thm.~13]{Fa03}, the topological complexity of $\b{S}^{\times n}$ is equal to $2n+1$.  Consequently, the motion planner constructed above is optimal.
\end{proof}

Now let $X$ be a subcomplex of the standard CW-decomposition of $\b{S}^{\times n}$.   In contrast to the odd-dimensional case, the topological complexity of $X$ is determined by the dimension of $X$.

\begin{thm} \label{thm:even subcomplex}
Let $\b{S}=S^{2k}$ be an even-dimensional sphere, and let $X$ be a $2k\ell$-dimensional subcomplex of the standard CW-decomposition of the product space $\b{S}^{\times n}$. Then the topological complexity of $X$ is $\tc(X)=2\ell +1$.
\end{thm}
\begin{proof} 
First we show that $X$ admits a motion planner with $2\ell+1$ local domains and rules. 
We assert that $(X \times X) \cap W_j=\emptyset$ for each $j$, $0 \le j \le 2n-2\ell-1$.  Given $\alpha
\in\set{0,1,2}^{\times n}$, identify $F_\alpha = F_{\alpha(1)} \times \cdots \times F_{\alpha(n)}$, 
and define $I_{\alpha}\subseteq [n]$ by $I_{\alpha}=\set{ i \in [n] \mid \alpha(i) <2}$.  That is, $I_\alpha$ is the set of indices $i$ for which $F_{\alpha(i)}$ contains only antipodal points in $\b{S}$, the $i$-th sphere in the product $\b{S}^{\times n}$.  

The cells of $X$ are of the form $A_J = \b{S}^{\times n}_J$ for certain subsets $J$ of $[n]$.  Consequently, $X\times X$ admits a CW-decomposition with cells $A_J \times A_K$ for certain subsets $J$ and $K$ of $[n]$.  If $\alpha$ satisfies $\abs{\alpha} < 2n-2\ell$, we can find $i \in I_{\alpha}$ with $i \notin J \cup K$.  In this situation, it is readily checked that $(A_J \times A_K)\cap F_{\alpha}=\emptyset$.   
Consequently, $\set{(X \times X) \cap W_j \mid 2n-2\ell \le j \le 2n}$ is a set of $2\ell+1$ local domains of a motion planner for $X$.

For $2n-2\ell \le j \le 2n$, let $\sigma_j(X)$ denote the restriction of the map $\sigma_j\colon W_j \to P\b{S}^{\times n}$ to $(X\times X) \cap W_j$.  Arguing as in the proof of Theorem \ref{thm:subcomplex}, one can show that that the image of $\sigma_j(X)$ is contained in $PX$.  It follows that the maps $\sigma_j(X)$, $2n-2\ell \le j \le 2n$, are local rules for the motion planner on $X$ with local domains given above.

Since the complex $X$ is $2k\ell$-dimensional, it must contain at least one copy of $\b{S}^{\times\ell}$ as a subcomplex.  Let $\tau\colon \b{S}^{\times\ell} \to X$ be the inclusion of such a subcomplex.  
Checking that this inclusion induces an epimorphism $\tau^* \colon H^*(X) \to H^*(\b{S}^{\times \ell})$ 
in cohomology, we conclude that $\zcl(H^*(X)) \ge \zcl(H^*(\b{S}^{\times\ell}))=2\ell$, and 
$\tc(X)\ge 2\ell +1$. Thus, the topological complexity of $X$ is equal to $2\ell +1$.
\end{proof}

\begin{remark} \label{rem:even sphere1}
If $X$ is $r$-connected, then 
\[
\tc(X) < \frac{2\cdot \dim(X)+1}{r+1} +1,
\]
see \cite[Thm.~5.2]{Fa04}.  This inequality may be used to verify that the 
topological complexity of a $2k\ell$-dimensional subcomplex $X$ of $\b{S}^{\times n}$ is bounded above by $2\ell+1$.  Note, however, that the proof of Theorem~\ref{thm:even subcomplex} provides an explicit optimal motion planner for $X$.
\end{remark}

\newcommand{\arxiv}[1]{{\texttt{\href{http://arxiv.org/abs/#1}{{#1}}}}}

\newcommand{\MRh}[1]{\href{http://www.ams.org/mathscinet-getitem?mr=#1}{MR#1}}

\bibliographystyle{amsplain}

\end{document}